\documentclass[12pt]{amsart}   
\usepackage{amscd,amssymb,amsmath}                  
\newtheorem{theorem}{Theorem}[section]
\newtheorem{lemma}[theorem]{Lemma}
\newtheorem{proposition}[theorem]{Proposition}
\newtheorem{question}[theorem]{Question}
\theoremstyle{definition}
\newtheorem{definition}[theorem]{Definition}
\theoremstyle{remark}

\numberwithin{equation}{section}

\begin{document} 
\noindent                                             
\begin{picture}(150,36)                               
\put(5,20){\tiny{Submitted to}}                       
\put(5,7){\textbf{Topology Proceedings}}              
\put(0,0){\framebox(140,34){}}                        
\put(2,2){\framebox(136,30){}}                        
\end{picture}                                         

\vspace{0.5in}

\title%
{On extremely amenable groups of homeomorphisms}

\author{Vladimir Uspenskij}

\address{321 Morton Hall, Department of Mathematics, Ohio
University, Athens, Ohio 45701, USA}

\email{uspensk@math.ohiou.edu}

\thanks{Research supported by BSF grant \#2006119}

\subjclass[2000]%
{Primary 54H11. Secondary 22A05, 22A15, 22F05, 37B05, 54H15, 54H20, 54H25, 57S05.}


\date{30 October 2007}

\newcommand{\abs}[1]{\lvert#1\rvert}
\def\norm#1{\left\Vert#1\right\Vert}
\def\Q {{\Bbb Q}}
\def\I {{\Bbb I}}
\def\C {{\Bbb C}}
\def\N{{\Bbb N}}
\def\di{{\mathrm{di}}}
\def\Z {{\Bbb Z}}
\def\U{{\Bbb U}}
\def\F{{\mathrm{E}}}
\def\Un{{\mathcal{U}}}
\def\Is{{\mathrm{Is}}\,}
\def\Aut{{\mathrm {Aut}}\,}
\def\supp{{\mathrm {supp}}\,}
\def\Homeo{{\mathrm{Homeo}}\,}
\def\gr{{\underline{\Box}}}
\def\diam{{\mathrm{diam}}\,}
\def\d{{\mathrm{dist}}}
\def\H{{\mathcal H}}
\def\me{{\mathrm{me}}}

\def\a{\alpha}
\def\d{\delta}
\def\D{\Delta}
\def\g{\gamma}
\def\s{\sigma}
\def\Si{\Sigma}
\def\implies{\Rightarrow}
\def\R{{\mathbf R}}
\def\Rp{{\mathbf R_+^*}}
\def\o{\omega}
\def\O{\Omega}
\def\G{\Gamma}

\def\sB{{\mathcal B}}
\def\sC{{\mathcal C}}
\def\sE{{\mathcal E}}
\def\sF{{\mathcal F}}
\def\sG{{\mathcal G}}
\def\sH{{\mathcal H}}
\def\sJ{{\mathcal J}}
\def\sK{{\mathcal K}}
\def\sM{{\mathcal M}}
\def\sN{{\mathcal N}}
\def\sO{{\mathcal O}}
\def\sS{{\mathcal S}}
\def\sT{{\mathcal T}}
\def\sU{{\mathcal U}}
\def\sV{{\mathcal V}}

\def\sbs{\subset}
\def\rar{\rightarrow}
\def\e{\epsilon}

\def\ti{\times}
\def\obr{^{-1}}
\def\stm{\setminus}
\def\newline{\hfill\break}

\def\Exp{{\mathrm{Exp}}\,}
\def\Iso{{\mathrm{Iso}}\,}
\def\Sym{{\mathrm{Sym}}\,}

\begin{abstract} 
A topological group $G$ is {\em extremely amenable} if every compact 
$G$-space has a $G$-fixed point. Let $X$ be compact and 
$G\sbs\Homeo(X)$.
We prove that the following are equivalent: 
(1)
$G$ is extremely amenable;
(2) 
every minimal closed
$G$-invariant subset of $\Exp R$ is a singleton, where $R$ is the closure
of the set of all graphs of $g\in G$ in the space $\Exp (X^2)$
($\Exp$ stands for the space of closed subsets);
(3) for each $n=1,2,\dots$ there is
a closed $G$-invariant subset $Y_n$ of $(\Exp X)^n$ 
such that 
$\cup_{n=1}^\infty Y_n$ contains arbitrarily fine covers of $X$
and
for every $n\ge1$ every minimal closed $G$-invariant subset of $\Exp Y_n$
is a singleton.
This yields
an alternative proof of Pestov's theorem that the group of all order-preserving 
self-homeomorphisms of the Cantor middle-third set (or of the interval $[0,1]$)
is extremely amenable.
\end{abstract}

\maketitle

\section{Introduction} 
\label{s:intro}

With every%
\footnote{All spaces are assumed to be Tikhonov,
and all {\em maps} are assumed to be continuous}  
topological group $G$ one can associate 
the {\em greatest ambit $\sS(G)$} and 
the {\em universal minimal compact $G$-space} $\sM(G)$. 
To define these objects, recall some 
definitions. A {\em $G$-space} is a topological space $X$ with a continuous
action of $G$, that is, a map $G\times X\to X$ satisfying $g(hx)=(gh)x$
and $1x=x$ ($g,h\in G$, $x\in X$). 
A map $f:X\to Y$ between two
$G$-spaces is {\em $G$-equivariant}, or a {\em $G$-map} for short, if
$f(gx)=gf(x)$ for every $g\in G$ and $x\in X$. 

A {\em semigroup} is a set with an associative multiplication.
A semigroup $X$ is {\em right topological} if it 
is a topological space and 
for every $y\in X$ the self-map $x\mapsto xy$ of $X$ is continuous. 
(Sometimes the term {\em left topological} is used for the same thing.)
A subset $I\sbs X$ is a {\em left ideal} if $XI\sbs I$. 
If $G$ is a topological group, a {\em right topological semigroup compactification}
of $G$ is a right topological compact semigroup $X$ together with a continuous 
semigroup morphism $f:G\to X$ with a dense range such that the map 
$(g,x)\mapsto f(g)x$ from $G\ti X$ to $X$ is jointly continuous (and hence 
$X$ is a $G$-space). 
The {\em greatest ambit} 
$\sS(G)$ for $G$ is a right topological semigroup compactification which is universal
in the usual sense: for any right topological semigroup compactification $X$ of $G$
there is a unique morphism $\sS(G)\to X$ of right topological semigroups
such that the obvious diagram commutes. Considered as a $G$-space, $\sS(G)$
is characterized by the following property:
there is a distinguished point $e\in \sS(G)$ such that for 
every compact $G$-space $Y$ and every $a\in Y$ there exists a unique 
$G$-map $f:\sS(G)\to Y$ such that $f(e)=a$. 

We can take for $\sS(G)$
the compactification of $G$ corresponding to the $C^\ast$-algebra $\hbox{RUCB}(G)$
of all bounded right uniformly continuous functions on $G$,
that is, the maximal ideal space of that algebra.
(A complex function $f$ on $G$ is {\it right uniformly continuous} if
$$
\forall\e>0\,\exists V\in\sN(G)\,\forall x,y\in G\,
(xy\obr\in V\implies|f(y)-f(x)|<\e),
$$
where $\sN(G)$ is the filter of neighbourhoods of unity.)
The $G$-space structure on $\sS(G)$ comes from the natural continuous action
of $G$ by automorphims on $\hbox{RUCB}(G)$ defined by $gf(h)=f(g\obr h)$ ($g,h\in G$, 
$f\in \hbox{RUCB}(G)$). We shall identify $G$ with
a subspace of $\sS(G)$. 
Closed 
$G$-subspaces of $\sS(G)$ are the same as closed left ideals of $\sS(G)$.

A $G$-space $X$ is {\em minimal} if it has
no proper $G$-invariant closed subsets or, equivalently, if the orbit
$Gx$ is dense in $X$ for every $x\in X$. 
The {\em universal minimal compact $G$-space} $\sM(G)$ is characterized
by the following property: $\sM(G)$ is a minimal compact $G$-space, and 
for every compact minimal $G$-space $X$ there exists a $G$-map of $\sM(G)$
onto $X$. Since Zorn's lemma implies that every compact $G$-space has 
a minimal compact $G$-subspace, it follows that for every compact $G$-space
$X$, minimal or not, there exist a $G$-map of $\sM(G)$ to $X$. The space $\sM(G)$
is unique up to a $G$-space isomorphism and is isomorphic to any minimal
closed left ideal of $\sS(G)$,
see e.g.~\cite{Aus}, \cite[Section 4.1]{Pestbook}, \cite[Appendix]{Usptopproc},
\cite[Theorem 3.5]{UspComp}.

A topological group $G$ is {\em extremely amenable}
if $\sM(G)$ is a singleton or, equivalently, if $G$ has the 
{\em fixed point on compacta property}: every compact $G$-space $X$
has a $G$-fixed point, that is, a point $p\in X$ such that $gp=p$ for
every $g\in G$. Examples of extremely amenable groups include
$\Homeo_+[0,1]=$ the group of all orientation-preserving self-homeomorphisms
of $[0,1]$; $U_s(H)=$ the unitary group of a Hilbert space $H$, with the topology
inherited from the product $H^H$; $\Iso(U)=$ the group of isometries of the 
Urysohn universal metric space $U$. See Pestov's book \cite{Pestbook}
for the proof. Note that a locally compact group $\ne\{1\}$
cannot be extremely amenable, since every locally compact group admits a free
action on a compact space \cite{Veech}, \cite[Theorem 3.3.2]{Pestbook}.

We refer the reader to Pestov's book \cite{Pestbook} for 
various intrinsic characterizations of extremely amenable groups.
These characterizations reveal 
a close connection between Ramsey theory and the notion of extreme
amenability. 
The aim of the present paper is to 
give another characterization of extremely  amenable groups,
based on a different approach. For a compact space $X$ let $H(X)$
be the group of all self-homeomorphisms of $X$, 
equipped with the compact-open topology. Let $G$ be a topological
subgroup of $H(X)$. There is an obvious necessary condition for $G$
to be extremely amenable: every minimal closed $G$-subset of $X$ must
be a singleton. However, this condition is not sufficient. For example,
let $X$ be the Hilbert cube, and let $G\sbs H(X)$ be the stabilizer of a given point
$p\in X$. Then the only minimal closed $G$-subset of $X$ is the singleton
$\{p\}$, but $G$ is not extremely amenable \cite{Usptopproc}, since $G$
acts without fixed points on the compact space $\Phi_p$ of all maximal
chains of closed subsets of $X$ starting at $p$. The space $\Phi_p$ is 
a subspace of the compact $G$-space $\Exp \Exp X$, where for a compact space
$K$ we denote by $\Exp K$ the compact space of all closed non-empty subsets
of $K$, equipped with the Vietoris topology%
\footnote{If $F$ is closed in $K$, the sets $\{A\in \Exp K: A\sbs F\}$ and 
$\{A\in \Exp K: A\hbox{ meets }F\}$ are closed in $\Exp K$, and the Vietoris
topology is generated by the closed sets of this form. If $K$ is a $G$-space,
then so is $\Exp K$, in an obvious way.}. 
It was indeed necessary to use
the second exponent in this example, the first exponent would not work.
One can ask whether in general for every group 
$G\sbs H(X)$ which is not extremely amenable there exists a compact $G$-space
$X'$ derived from $X$ by applying a small number of simple functors, like
powers, probability measures, exponents, etc., such that $X'$ contains a
closed $G$-subspace (which can be taken minimal) on which $G$ acts without
fixed points. We answer this question in the affirmative. 

Consider the action of $G$ on $\Exp(X^2)$ defined by the composition
of relations: if $g\in G$, $F\sbs X^2$, and $\G_g\sbs X^2$ is the graph of $g$, 
then $gF=\G_g\circ F=\{(x,gy):(x,y)\in F\}$. This amounts to considering $X^2$
as the product of two different $G$-spaces: the first copy of $X$ has 
the trivial $G$-structure, and the second copy is the given $G$-space $X$.
If $G$ is not extremely amenable, then there is a closed minimal
$G$-subspace $Y$ of $\Exp\Exp (X^2)$ that is not a singleton (and hence
fixed point free). This follows from:

\begin{theorem}
\label{th1}
Let $X$ be compact, $G$ a subgroup of $H(X)$. Denote by $R$
the closure of the set $\{\G_g:g\in G\}$ of the graphs
of all $g\in G$ in the space $\Exp (X^2)$.
Then $G$ is extremely amenable if and only if every
minimal closed $G$-subset of $\Exp R$ is a singleton. 
\end{theorem}

Here $X^2$ is the product of the trivial $G$-space and the given $G$-space $X$,
as in the paragraph preceding Theorem~\ref{th1}, and $R$ is considered as a $G$-subspace
of $\Exp (X^2)$. 

For example, let $X=I=[0,1]$ be the closed unit interval. Consider
the group $G=H_+([0,1])$ of all orientation-preserving self-homeomorphisms of $I$.
The space $R$ in this case consists of all curves $\G$ in the square $I^2$ that
connect the lower left and upper right corners and ``never go down": 
if $(x,y)\in \G$, $(x',y')\in \G$ and $x<x'$, then $y\le y'$
(see the picture in \cite[Example 2.5.4]{PWhere}). It can be verified
that the only minimal compact $G$-subsets of $\Exp R$ are singletons 
(they are of the form \{a closed union of $G$-orbits in $R$\}). The proof depends on the following lemma:

\begin{lemma}
\label{l1}
Let $\D^n$ be the $n$-simplex of all $n$-tuples $(x_1,\dots, x_n)\in I^n$ such that
$0\le x_1\le \dots \le x_n\le 1$. Equip $\D^n$ with the natural action of
the group $G=H_+([0,1])$. Then every minimal closed $G$-subset of $\Exp \D^n$
is a singleton (= \{a union of some faces of $\D^n$\}).
\end{lemma}

The idea to consider the action of $G=H_+([0,1])$ on $\D^n$ is borrowed from
\cite{Drin}, where it is shown that the geometric realization of any simplicial
set can be equipped with a natural action of $G$. We shall not prove Lemma~\ref{l1},
since this lemma follows from Pestov's theorem that $G$ is extremely amenable, and 
I am not aware of a short independent proof of the lemma. The essence of the lemma
is that every subset of $\D^n$ can be either pushed (by an element of $G$) into
the $\e$-neighbourhood of the boundary of the simplex or else can be pushed 
to approximate the entire simplex within $\e$. Some Ramsey-type argument
seems to be necessary for this. Actually 
Lemma~\ref{l1}
may be viewed as a topological equivalent of the finite Ramsey theorem
\cite[Theorem 1.5.2]{Pestbook}, since 
Pestov showed that this theorem has an equivalent reformulation in terms
of the notion of a ``finitely oscillation stable" dynamical system 
\cite[Section 1.5]{Pestbook}, and extremely amenable groups are characterized
in the same terms \cite[Theorem 2.1.11]{Pestbook}.

An important example of an extremely amenable group is the Polish group $\Aut(\Q)$ of all
automorphisms of the ordered set $\Q$ of rationals \cite{P-Ellis},
\cite[Theorem 2.3.1]{Pestbook}. This group is considered with the topology inherited from 
$(\Q_d)^\Q$, where $\Q_d$ is the set of rationals with the discrete topology. Let $K\sbs [0,1]$
be the usual middle-third Cantor  set. The topological group $\Aut(\Q)$ is isomorphic to the 
topological group $G=H_<(K)\sbs H(K)$ of all order-preserving self-homeomorphisms of $K$.
To see this, note that pairs of the endpoints of  ``deleted intervals" (= components of $[0,1]\stm K$)
form a set which is order-isomorphic to $\Q$, whence a homomorphism $G\to \Aut(\Q)$ which
is easily verified to be a topological isomorphism. 
One can prove that  the group $G\simeq\Aut(\Q)$
is extremely amenable with the aid of Theorem~\ref{th1}. The proof is essentially the same as 
in the case of the group $G=H_+([0,1])$. 
The space $R$ considered in Theorem~\ref{th1} again is the space of  ``curves", this time
in $K^2$, that go from $(0,0)$ to $(1,1)$ and ``look like graphs", with the exception that
they may contain vertical and horizontal parts. The evident analogue of Lemma~\ref{l1}
holds for ``Cantor simplices" of the form 
$\{(x_1,\dots, x_n)\in K^n:0\le x_1\le \dots \le x_n\le 1\}$.

Theorem~\ref{th1} may help to answer the following:

\begin{question}
\label{q1}
Let $P$ be pseudoarc, $G=H(P)$, and let $G_0$ be the stabilizer of a given
point $x\in P$. Is $G_0$ extremely amenable? 
\end{question}

As explained in \cite{Usptopproc}, this question is motivated by the
observation that the argument involving maximal chains, which shows
that the stabilizer $G_0\sbs H(X)$ of a point $p\in X$
is not extremely amenable if $X$ is either a Hilbert cube 
or a compact manifold of dimension $>1$, does not
work for the pseudoarc. A positive answer to Question~\ref{q1} would
imply that the pseudoarc $P$ can be identified with $\sM(G)$ for $G=H(P)$.
The problem whether this is the case was raised in \cite{Usptopproc}
and appears as Problem 6.7.20 in \cite{Pestbook}.

The {\em suspension} $\Si X$ 
of a space $X$ is the quotient of $X\ti I$ obtained
by collapsing the ``bottom" $X\ti\{0\}$ and the ``top" $X\ti \{1\}$ to points.
Let $q:\Si X\to I$ be the natural projection. The inverse image under $q$ of
the maximal chain $\{[0,x]:x\in I\}$ of closed subsets of $I$ 
is a maximal chain of closed subsets of $\Si X$.

\begin{question}
\label{q2}
Let $Q=I^\omega$ be the Hilbert cube, and $C$ be the 
maximal chain of subcontinua of $\Si Q$ considered above. If $G=H(\Si Q)$
and $G_0\sbs G$ is the stabilizer of $C$, is $G_0$ extremely amenable?
\end{question}

This question is motivated by the search for a good candidate for the
space $\sM(G)$, where $G=H(Q)$. The space $\Phi_c$ of all maximal chains of
subcontinua of $Q$, proved to be minimal by Y. Gutman \cite{Gutman}, 
may be such a candidate \cite[Problem 6.4.13]{Pestbook}. 
Recall that for the group $G=H(K)$, where $K=2^\o$
is the Cantor set, $\sM(G)$ can be identified with the space $\Phi\sbs \Exp\Exp K$
of all maximal chains of closed subsets of $K$ \cite{GW}, 
\cite[Example 6.7.18]{Pestbook}. 

There is another characterization 
(Theorem~\ref{th:new})
of extremely amenable groups in the spirit
of Theorem~\ref{th1} which, in combination with Lemma~\ref{l1}, readily implies
Pestov's results that $H_+([0,1])$ 
and $\Aut(\Q)$ are extremely amenable. Let $X$ be compact,
$Y_n\sbs (\Exp X)^n$ for $n=1,2,\dots$. We say that $\cup_{n=1}^\infty Y_n$ 
{\em contains arbitrarily fine covers} if for every open cover $\a$ of $X$
there are $n\ge1$ and $(F_1, \dots, F_n)\in Y_n$ such that $\cup_{i=1}^n F_i=X$
and the cover $\{F_i\}_{i=1}^n$ of $X$ refines $\a$. 

\begin{theorem}
\label{th:new}
Let $X$ be compact, $G$ a subgroup of $H(X)$. Let
$Y_n$ be a closed $G$-invariant subset of $(\Exp X)^n$ ($n=1,2,\dots$)
such that $\cup_{n=1}^\infty Y_n$ contains arbitrarily fine covers of $X$. 
Then $G$ is extremely amenable if and only if
for every $n\ge1$ every minimal closed $G$-invariant subset of $\Exp Y_n$
is a singleton.
\end{theorem}

Observe that Pestov's theorem asserting that $G=H_+([0,1])$ is extremely
amenable follows from Theorem~\ref{th:new} and Lemma~\ref{l1}:
it suffices to take for $Y_{n+1}$ the collection of all sequences
$$
([0,x_1], [x_1, x_2], \dots, [x_n, 1]),
$$
where 
$0\le x_1\le \dots \le x_n\le 1$. The $G$-space $Y_{n+1}$ is isomorphic
to the $n$-simplex $\D^n$ considered in Lemma~\ref{l1}.
The argument for $\Aut(\Q)\simeq H_<(K)$ is similar.

The proof of Theorems~\ref{th1} 
and~\ref{th:new} depends on the notion of a representative
family of compact $G$-spaces. We introduce this notion in Section~\ref{s:2}
and observe that a topological group $G$ is extremely amenable if (and only if) 
there exists a representative family $\{X_\a\}$ such that any minimal closed 
$G$-subset of any $X_\a$ is a singleton (Theorem~\ref{th:crit}). 
In Section~\ref{s:3} we prove that the single space $\Exp R$ considered in 
Theorem~\ref{th1} constitutes a representative family (Theorem~\ref{th:x2}).
The conjunction of Theorems~\ref{th:crit} and~\ref{th:x2} proves Theorem~\ref{th1}. 
In Section~\ref{s:4} we prove that under the conditions of Theorem~\ref{th:new}
the sequence $\{\Exp Y_n\}$ is representative (Theorem~\ref{th:another}).
The conjunction of Theorems~\ref{th:crit} and~\ref{th:another} proves 
Theorem~\ref{th:new}.

\section{Representative families of $G$-spaces}
\label{s:2}

Let $G$ be a topological group, $X$ a compact $G$-space. 
For $g\in G$ the {\em $g$-translation} of $X$
is the map $x\mapsto gx$,  $x\in X$. 
The {\em enveloping semigroup} (or the {\em Ellis semigroup}) 
$E(X)$ of the dynamical system $(G,X)$ is the closure
of the set of all $g$-translations, $g\in G$, in the compact space $X^X$.
This is a right topological semigroup compactification of $G$, as defined in Section~\ref{s:intro}.
The natural
map $G\to E(X)$ extends to a $G$-map $\sS(G)\to E(X)$ which is a morphism
of right topological semigroups. 

\begin{definition}
A family $\{X_\a:\a\in A\}$ of compact $G$-spaces is {\em representative}
if the family of natural maps $\sS(G)\to E(X_\a)$, $\a\in A$, separates
points of $\sS(G)$ (and hence yields an embedding of $\sS(G)$ into 
$\prod_{\a\in A} E(X_\a)$).
\end{definition}

\begin{theorem}
\label{th:crit}
Let $G$ be a topological group, $\{X_\a\}$ 
a representative family of compact $G$-spaces.
Then $G$ is extremely amenable if (and only if) 
every minimal closed $G$-subset of every $X_\a$ is a singleton. 
\end{theorem}

This is a special case of a more general theorem:

\begin{theorem}
\label{th:gen}
If $\{X_\a\}$  is a representative family of compact $G$-spaces,
the universal minimal compact $G$-space $\sM(G)$ is isomorphic (as a $G$-space)
to a $G$-subspace of a product $\prod Y_\beta$, where each $Y_\beta$ is a minimal
compact $G$-space isomorphic to a $G$-subspace of some $X_\a$.
\end{theorem}

\begin{proof}
By definition of a representative family, the greatest ambit $\sS(G)$ can be embedded
(as a $G$-space) into the product $\prod E(X_\a)$ and hence also into the product
$\prod X^{X_\a}_\a$. Consider $\sM(G)$ as a subspace of $\sS(G)$ and take for 
the $Y_\beta$ 's the projections of $\sM(G)$ to the factors $X_\a$. 
\end{proof}

We now give a sufficient condition for a family of compact $G$-spaces to be representative.
Let us say that two subsets $A,B$ of $G$ are 
{\em far from each other with respect to the right uniformity} if one of the following 
equivalent conditions holds: (1) the neutral element $1_G$ of $G$ 
is not in the closure of the set $BA\obr$; (2)  
for some neighbourhood $U$ of $1_G$ the sets $A$ and $UB$ are disjoint;
(3) there exists a right uniformly continuous function $f: G\to [0,1]$ such that
$f=0$ on $A$ and $f=1$ on $B$; (4) $A$ and $B$ have disjoint closures in $\sS(G)$.

\begin{proposition}
\label{p1}
Let $\sF$ be a family of compact $G$-spaces. Suppose that the following holds:

(*) if $A,B\sbs G$ are far from each other with respect to the right uniformity,
then there exists $X\in \sF$ and $p\in X$ such that the sets $Ap$ and $Bp$ have
disjoint closures in $X$.

Then $\sF$ is representative.
\end{proposition}

\begin{proof}
Consider the natural map $G\to \prod\{E(X):X\in \sF\}$. 
It defines a compactification $bG$ of $G$. We must prove that this compactification 
is equivalent to $\sS(G)$. 

Let $A,B$ be any two subsets  of $G$ with disjoint closures in $\sS(G)$.
Then $A$ and $B$ are far from each other 
with respect to the right uniformity. According to the condition (*), 
there exists $X\in \sF$ and $p\in X$ such that the sets $Ap$ and $Bp$ have
disjoint closures in $X$. It follows that the images of $A$ and $B$ in $E(X)$
have disjoint closures, and {\em a fortiori} the images of $A$ and $B$ in $bG$
have disjoint closures. It follows that $\sS(G)$ and $bG$ are equivalent compactifications 
of $G$ \cite[Theorem 3.5.5]{E}.
\end{proof}

\section{Proof of Theorem~\ref{th1}}
\label{s:3}

Recall the setting of Theorem~\ref{th1}: $X$ is compact, $G$ is a topological subgroup of 
$H(X)$. For $g\in G$ let  $\G_g=\{(x,gx):x\in X\}\sbs X^2$ be the graph of $g$, and let
 $R$  be the closure of the set $\{\G_g:g\in G\}$ in the compact space $\Exp(X^2)$. 
We consider the action of $G$ on $\Exp(X^2)$ defined by $gF=\{(x,gy):(x,y)\in F\}$ 
($g\in G$, $F\in \Exp(X^2)$), and consider $R$ as a $G$-subspace of $\Exp(X^2)$.

\begin{theorem}
\label{th:x2}
Let $X$ be a compact space, $G\subset H(X)$.
Let $R\sbs \Exp(X^2)$ be the compact $G$-space defined above.
The family consisting of the single compact $G$-space 
$\Exp R$ 
is representative.
\end{theorem}

In other words, $\sS(G)$ is isomorphic to the enveloping semigroup of $\Exp R$.

\begin{proof}
Let $A, B\sbs G$ be far from each other (that is, $1_G$ is not in the closure
of $BA\obr$). In virtue of proposition~\ref{p1}, it suffices to find 
$p\in Y=\Exp R$ such that $Ap$ and $Bp$ have disjoint 
closures in $Y$. 

Let $p$ be the closure of the set $\{\G_g:g\in A\obr\}$ in the space 
$\Exp(X^2)$. Then $p$ is a closed subset of $R$ and hence $p\in Y$.
We claim that $p$ has the required property: $Ap$ and $Bp$ have disjoint
closures in $Y$ or, which is the same, in $\Exp\Exp (X^2)$.

There exist a continuous pseudometric $d$ on $X$ and $\d>0$ such that
$$
\forall f\in A\ \forall g\in B\ \exists x\in X \ (d(gf\obr(x), x)\ge \d).
$$
Let $\D\sbs X^2$ be the diagonal. Let $C\sbs X^2$ be the closed set 
defined by 
$$
C=\{(x,y)\in X^2: d(x,y)\ge \d\}.
$$
Let $K\sbs \Exp X^2$ be the closed
set defined by 
$$
K=\{F\sbs X^2: F\hbox{ meets }C\}.
$$
Consider the closed
sets $L_1,L_2\sbs \Exp\Exp (X^2)$ defined by 
$$
L_1=\{q\sbs \Exp(X^2):q\hbox{ is closed and }\D\in q\}
$$
and
$$
L_2=\{q\sbs \Exp(X^2):q\hbox{ is closed and }q\sbs K\}.
$$
Since $\D\notin K$, the sets $L_1$ and $L_2$ are disjoint. 
It suffices to 
verify that $Ap\sbs L_1$ and $Bp\sbs L_2$.

The first inclusion is immediate: if $g\in A$, then for $h=g\obr$
we have $\D=g\G_h\in gp$, hence $gp\in L_1$. Thus $Ap\sbs L_1$.
We now prove that $Bp\sbs L_2$. Let $g\in B$. If $f\in A$ and $h=f\obr$,
there exists $x\in X$ such that $d(gh(x), x)\ge \d$, which means that 
$\G_{gh}$ meets $C$. Hence $g\G_h=\G_{gh}\in K$. It follows that the closed
set $g\obr K$ contains the set $\{\G_h: h\in A\obr\}$ and hence also its closure
$p$. In other words, $gp\sbs K$ and hence $gp\in L_2$.
\end{proof}

As noted in Section~\ref{s:intro}, Theorem~\ref{th1} follows from 
Theorems~\ref{th:crit} and~\ref{th:x2}.

Combining Theorems~\ref{th:gen} and~\ref{th:x2}, we obtain the following
generalization of Theorem~\ref{th1}:

\begin{theorem}
\label{c3}
Let $X$ be a compact space, $G$ a subgroup of $H(X)$.
Let $R$ be the same as in Theorems~\ref{th1} and~\ref{th:x2}.
Let $\sF$ be the family of all 
minimal closed $G$-subspaces of $\Exp R$.
Then $\sM(G)$ is isomorphic to a subspace of a product
of members of $\sF$ (some factors may be repeated).
\end{theorem}

\section{Proof of Theorem~\ref{th:new}}
\label{s:4}

Theorem~\ref{th:x2} implies that for any subgroup $G\sbs H(X)$
the one-point family $\{\Exp\Exp(X^2)\}$ is representative
(recall that we consider the trivial action on the first factor $X$).
I do not know whether $X^2$ can be replaced here by $X$.
On the other hand, the following holds:

\begin{theorem}
\label{th:oldanother}
Let $X$ be a compact space, $G$ a subgroup of $H(X)$.
The sequence \newline
$\{\Exp((\Exp X)^n)\}_{n=1}^\infty$ of compact $G$-spaces
is representative.
\end{theorem}

This is a special case of a more general theorem:

\begin{theorem}
\label{th:another}
Let $X$ be a compact space, $G$ a subgroup of $H(X)$.
Let
$Y_n$ be a closed $G$-invariant subset of $(\Exp X)^n$ ($n=1,2,\dots$)
such that $\cup_{n=1}^\infty Y_n$ contains arbitrarily fine covers of $X$.
Then the sequence 
$\{\Exp Y_n\}_{n=1}^\infty$ of compact $G$-spaces
is representative.
\end{theorem}

\begin{proof}
Let $A,B\sbs G$ be two sets that are far from each other with respect
to the right uniformity. In virtue of proposition~\ref{p1}, 
it suffices to find $n$ and a point 
$p\in \Exp Y_n$ such that $Ap$ and $Bp$ have disjoint 
closures in $\Exp Y_n$ or, which is the same, in $Z_n=\Exp((\Exp X)^n)$.

There exist a continuous pseudometric $d$ on $X$ and $\d>0$ such that
$A$ and $B$ are $(d,2\d)$-far from each other, in the sense that
$$
\forall f\in A\ \forall g\in B\ \exists x\in X \ (d(f(x), g(x))>2\d).
$$
The assumption that $\cup_{n=1}^\infty Y_n$ contains arbitrarily fine covers
implies that we can find $n\ge 1$ and closed sets 
$C_1, \dots, C_n\sbs X$ of $d$-diameter
$\le \d$ such that $(C_1, \dots, C_n)\in Y_n$ and $\cup_{i=1}^n C_i=X$.
For each $g\in G$ let 
$F_g=(g\obr(C_1),\dots,g\obr(C_n))\in (\Exp X)^n$. 
Since $Y_n$ is $G$-invariant, we have $F_g\in Y_n$. 
Let $p$ be the closure
of the set $\{F_g:g\in A\}$ in the space $(\Exp X)^n$. Then $p\in \Exp Y_n$.
We claim that $p$ has the required property: $Ap$ and $Bp$ have disjoint
closures in $Z_n$. 

Let $D_i=\{x\in X:d(x,C_i)\ge\d\}$, $i=1,\dots,n$.
Consider the closed sets $K_1, K_2\sbs (\Exp X)^n$ defined by 
$$
K_1=\{(F_1,\dots,F_n)\in (\Exp X)^n: F_i\sbs C_i,\ i=1,\dots,n\}
$$
and 
$$K_2=\{(F_1,\dots,F_n)\in (\Exp X)^n: F_i \hbox{ meets }D_i 
\hbox{ for some } i=1, \dots, n\}.
$$
Consider the closed sets $L_1, L_2\sbs Z_n$ defined by 
$$
L_1=\{q\sbs (\Exp X)^n:q\hbox{ is closed and }q\hbox{ meets }K_1\}
$$ 
and 
$$
L_2=\{q\sbs (\Exp X)^n:q\hbox{ is closed and }q\sbs K_2\}.
$$
Clearly $K_1$ and $K_2$ are
disjoint, hence $L_1$ and $L_2$ are disjoint as well. It suffices to 
verify that $Ap\sbs L_1$ and $Bp\sbs L_2$.

The first inclusion is immediate: if $g\in A$, then $F_g\in p$ and
$gF_g=(C_1,\dots,C_n)\in K_1\cap gp$, hence $gp$ meets $K_1$ and
$gp\in L_1$. We now prove that 
$Bp\sbs L_2$. Let $h\in B$. If $g\in A$, we can find $x\in X$
such that $d(g(x), h(x))>2\d$ and an index $i$, $1\le i\le n$, such
that $g(x)\in C_i$. Since $\diam C_i\le \d$, we have $h(x)\in D_i$ 
and therefore 
$h(x)\in hg\obr(C_i)\cap D_i\ne\emptyset$. 
It follows that 
$hF_g=(hg\obr(C_1),\dots,hg\obr(C_n))\in K_2$. This holds for every
$g\in A$, and thus we have shown that the closed set 
$h\obr K_2\sbs(\Exp X)^n$ contains the set $\{F_g:g\in A\}$ and
hence also its closure $p$. In other words, $hp\sbs K_2$ and hence
$hp\in L_2$.
\end{proof}

Theorem~\ref{th:new} follows from Theorems~\ref{th:another}
and~\ref{th:crit}.

Combining Theorems~\ref{th:gen} and~\ref{th:another}, we obtain
the following generalization of Theorem~\ref{th:new}:

\begin{theorem}
\label{c1}
Let $X$ be a compact space, $G$ a subgroup of $H(X)$.
Let
$Y_n$ be a closed $G$-invariant subset of $(\Exp X)^n$ ($n=1,2,\dots$)
such that $\cup_{n=1}^\infty Y_n$ contains arbitrarily fine covers of $X$.
Let $\sF$ be the family of all (up to an isomorphism)
minimal closed $G$-subspaces of $\Exp Y_n$, $n=1,2,\dots$.
Then $\sM(G)$ is isomorphic to a subspace of a product
of members of $\sF$ (some factors may be repeated).
\end{theorem}

\bibliographystyle{plain}

\end{document}